\documentclass[reqno,11pt]{amsart}
\usepackage{psfig, amsmath, amsfonts, amssymb, amsthm, amscd}

\DeclareMathSymbol{\leqslant}{\mathalpha}{AMSa}{"36} 
\DeclareMathSymbol{\geqslant}{\mathalpha}{AMSa}{"3E} 
\DeclareMathSymbol{\eset}{\mathalpha}{AMSb}{"3F}     
\renewcommand{\leq}{\;\leqslant\;}                   
\renewcommand{\geq}{\;\geqslant\;}                   
\newcommand{\dd}{\,\text{\rm d}}             
\newcommand{\sumtwo}[2]{\sum_{\substack{#1 \\ #2}}} 
\newcommand{\limtwo}[2]{\lim_{\substack{#1 \\ #2}}}     

\makeatletter
\def\captionfont@{\footnotesize}
\def\captionheadfont@{\scshape}

\long\def\@makecaption#1#2{%
  \vspace{2mm}
  \setbox\@tempboxa\vbox{\color@setgroup
    \advance\hsize-6pc\noindent
    \captionfont@\captionheadfont@#1\@xp\@ifnotempty\@xp
        {\@cdr#2\@nil}{.\captionfont@\upshape\enspace#2}%
    \unskip\kern-6pc\par
    \global\setbox\@ne\lastbox\color@endgroup}%
  \ifhbox\@ne 
    \setbox\@ne\hbox{\unhbox\@ne\unskip\unskip\unpenalty\unkern}%
  \fi
  \ifdim\wd\@tempboxa=\z@ 
    \setbox\@ne\hbox to\columnwidth{\hss\kern-6pc\box\@ne\hss}%
  \else 
    \setbox\@ne\vbox{\unvbox\@tempboxa\parskip\z@skip
        \noindent\unhbox\@ne\advance\hsize-6pc\par}%
\fi
  \ifnum\@tempcnta<64 
    \addvspace\abovecaptionskip
    \moveright 3pc\box\@ne
  \else 
    \moveright 3pc\box\@ne
    \nobreak
    \vskip\belowcaptionskip
  \fi
\relax
}
\makeatother
\def\writefig#1 #2 #3 {\rlap{\kern #1 truecm
\raise #2 truecm \hbox{#3}}}

\psfigurepath{.:./pictures}

\newtheorem{thm}{Theorem}[section]

\newcommand{\cB}{\ensuremath{\mathcal B}}
\newcommand{\cC}{\ensuremath{\mathcal C}}
\newcommand{\cD}{\ensuremath{\mathcal D}}

\newcommand{\cX}{\ensuremath{\mathcal X}}

\newcommand{\frS}{\ensuremath{\mathfrak S}}
\newcommand{\frT}{\ensuremath{\mathfrak T}}

\newcommand{\bbR}{{\ensuremath{\mathbb R}} }

\newcommand{\bbZ}{{\ensuremath{\mathbb Z}} }

\newcommand{\gd}{\delta}
\newcommand{\gep}{\varepsilon}       

\newcommand{\go}{\omega}

\newcommand{\gL}{\Lambda}

\newcommand{\setof}[2]{{\{#1\,:\,#2\}}}
\newcommand{\bnd}{\partial}

\newcommand{\bk}[1]{{\langle#1\rangle}}
\newcommand{\abs}[1]{\lvert#1\rvert}

\newcommand{\IF}[1]{1_{\{#1\}}}

\newcommand{\thebox}{{\gL_N}}
\newcommand{\Ac}{{\thebox\setminus A}}

\newcommand{\nn}{\sim}

\begin{document}
\title{A note on wetting transition for gradient fields}
\author{P.~Caputo and Y.~Velenik}
\address{
Fachbereich Mathematik, Sekr. MA 7-4, TU-Berlin, Stra\ss e des 17. Juni 136,
D-10623 Berlin, Germany
}
\email{caputo\@@math.tu-berlin.de}
\email{velenik\@@math.tu-berlin.de}
\date{\today}
\begin{abstract}
We prove existence of a wetting transition for two types of gradient fields: 1)
Continuous SOS models in any dimension and 2) Massless Gaussian model in
dimension $2$. Combined with a recent result showing the absence of such a
transition for Gaussian models above $2$ dimensions \cite{BDZ}, this shows in
particular that absolute-value and quadratic interactions can give rise to
completely different behavior.
\end{abstract}
\maketitle
\section{Introduction}
In several recent papers \cite{BoIo, DeVe, BoBr, BoStFlour, BMF}, the question has
been raised whether the 2 dimensional massless Gaussian model exhibited a
wetting transition. It is well-known that the Gaussian field in 2D is
delocalized, with a logarithmically divergent mean height, but that the
introduction of an arbitrarily weak self-potential favoring height $0$ is
enough to localize it, in the sense that the mean height remains finite
\cite{DMRR,BoBr}; this result has recently been extended to a class of non
Gaussian models in a stronger form, showing in particular existence of
exponential moments for the heights \cite{DeVe} and exponential decay of
covariances \cite{IoVe}. For higher dimensional model, the field is already
localized without pinning potential, but the introduction of such a potential turns
the algebraic decay of the covariances into an exponential one. On the other
hand, a Gaussian field with a positivity constraint (``surface above a
hard-wall'') exhibits entropic repulsion: The average height diverges like
$\log N$ in 2D \cite{2D} and 
$\sqrt{\log N}$ in higher dimensions \cite{BDZ1} ($N$
being the linear size of the box). When both a positivity constraint and a
pinning potential are present (``surface above an attractive hard wall''), 
there
is a competition between these two effects. If there exists a (non zero)
critical value for the strength of the pinning potential above which the
interface is localized, but below which it is repelled by the wall, we say that
the model exhibits a wetting transition. That such a transition occurs in a wide
class of 1D model is well-known, see e.g. \cite{Fi,BDZ}. It was recently shown
in~\cite{BDZ} that the Gaussian model in dimensions $3$ or higher does not
display a wetting transition: The interface is always localized. The physically
important case of the 2D model (describing a 2D interface in a 3D
medium) remained however open.

In the present note, we prove that the 2D Gaussian field does exhibit a
wetting transition; in fact the proof applies to any strictly convex
interaction, see below. We also prove that the continuous SOS model has such a
transition in {\em any} dimension, thus showing that the choice of the
interaction can greatly affect the physics of the system. Our proofs are based
on a variant of a beautiful and simple argument of Chalker 
\cite{Chalker}, who proved
the existence of a wetting transition in the discrete SOS model in dimension 2.
\section{Results}
We consider a class of gradient models with single spin-space $\bbR^+$, i.e.
modeling surfaces above a hard wall. Let
$\thebox$ be the cube of side $N$ centered at the origin, 
and $\Psi:\bbR\to\bbR$ an even
function to be specified later; we consider the following Hamiltonian
$$
H_N^{0,a,b}(\phi) = H_{0,N}^0 (\phi) + V_N^{a,b}(\phi)\,,
$$
where
\begin{align*}
H_{0,N}^0 (\phi) &= \sum_{\bk{x,y}\subset\thebox} \Psi(\phi_x-\phi_y)
\;+\sumtwo{\bk{x,y}}{x\in\thebox,\,y\notin\thebox} \Psi(\phi_x)\,,\\
V_N^{a,b}(\phi) &= -b\sum_{x\in\thebox} \IF{\phi_x\leq a}\,, \;\;\;a,b > 0\,,
\end{align*}
($\bk{x,y}$ denotes a pair of nearest neighbour sites).
The corresponding Gibbs measure (on $(\bbR^+)^\thebox$) is then given by
\begin{equation*}
\mu_N^{0,+,a,b}(\dd\phi) = \frac{e^{-H_N^{0,a,b}(\phi)}}{Z_N^{0,+,a,b}}
\prod_{x\in\thebox} \dd\phi_x\,.
\end{equation*}
As in the pure pinning problem (i.e. without a wall) 
\cite{DeVe}, the
relevant parameter is $\gep(a,b)=ae^b$ and not both $a$ and $b$ separately.
As usual, we introduce the {\em $\gd$-pinning} limit, which is the model
described by the measure
$$
\mu_N^{0,+,\gep}(\dd\phi) = \limtwo{a\to 0}{\gep(a,b)=\gep} 
\mu_N^{0,+,a,b}(\dd\phi) = \frac{e^{-H_{0,N}^0(\phi)}}{Z_N^{0,+,\gep}}
\prod_{x\in\thebox} (\dd\phi_x + \gep \gd_0(\dd\phi_x))\,.
$$
Note that $\mu_N^{0,+,\gep}$ can be written more explicitly as
$$
\mu_N^{0,+,\gep}(\dd\phi) = \sum_{A\subset\thebox} \gep^{\abs A}
\frac{Z^{0,+}_\Ac}{Z_N^{0,+,\gep}}\,\mu^{0,+}_{\Ac}(\dd\phi)\,,
$$
where
\begin{align*}
Z^{0,+}_\Ac = \int e^{-H_{0,N}^0(\phi)}\,\prod_{x\in\Ac} \dd\phi_x \prod_{y\in
A} \gd_0(\dd\phi_y)\,,\\
\intertext{and}
\mu^{0,+}_{\Ac}(\dd\phi) = (Z^{0,+}_\Ac)^{-1}\, e^{-H_{0,N}^0(\phi)}
\prod_{x\in\Ac} \dd\phi_x \prod_{y\in A} \gd_0(\dd\phi_y)
\end{align*}
is the Gibbs measure on $\Ac$ with $0$ boundary condition outside.
%
\smallskip

\noindent
{\bf Remark:} Here and everywhere else in this note, the integrals are
restricted to the positive real axis, so we do not write this condition
explicitly.

\smallskip
A quantity of interest is the density of {\em pinned sites}, i.e. of those sites
where the interface feels the effect of the pinning potential; it is defined as
$$
\rho_N(a,b) = \abs\thebox^{-1} \mu_N^{0,+,a,b} (\nu_N(\phi))\,,
$$
where $\nu_N(\phi) = \sum_{x\in\thebox}\IF{\phi_x\leq a}$. We also write
$\rho(a,b)=\lim_{N\to\infty} \rho_N(a,b)$. The corresponding quantities in the
$\gd$-pinning limit are denoted $\rho_N(\gep)$, $\rho(\gep)$
(measuring the density of sites exactly at height $0$). 
$\rho$ will play the role of an order parameter
for the wetting transition. 

Our results can then be stated as follows.
\begin{thm}\label{thm_abs}
Let $\Psi(x)=\abs x$, $d \geq 1$. 
Suppose that $ae^b < (2d)^{-1}$, then there exist two
constants $C_1(a,b,d)$ and $C_2(a,b,d)$ such that, for any $M>C_1\,N^{d-1}$,
$$
\mu_N^{0,+,a,b} (\nu_N(\phi)>M) \leq e^{-C_2\,M}\,.
$$
In particular, $\rho(a,b)=0$.
\end{thm}
%
%
\begin{thm}\label{thm_deltapinning}
Consider $\gd$-pinning in the two cases:
\begin{enumerate}
\item{
$\Psi(x)= \abs x$, $d \geq 1$, $\gep < e^{-2d}$;
}
\item{$\Psi(x)=\frac12 x^2$, $d=2$, $\gep < \gep_0$ for
some small positive number $\gep_0$.
}
\end{enumerate}
Then there exist two
constants $C_3$ and $C_4$ 
depending on $\gep,d$ in the first case and only 
on $\gep$ in the second case,
such that, for any $M>C_3\,N^{d-1}$,
$$
\mu_N^{0,+,\gep} (\nu_N(\phi)>M) \leq e^{-C_4\,M}\,.
$$
In particular, $\rho(\gep)=0$ in both cases.
\end{thm}
We recall that it is not hard
to prove that $\rho>0$ when the pinning is strong. For example, in
the $\gd$-pinning case, we can proceed in the following way. Since
$$
\abs\thebox^{-1}\log\frac{Z_N^{0,+,\gep}}{Z_N^{0,+,0}} = \int_0^\gep
\frac1{\hat\gep}\,\rho_N(\hat\gep)\,\dd \hat\gep\,,
$$
the result follows from $Z_N^{0,+,\gep}\geq\gep^{\abs\thebox}$ and the
existence of a constant $C$ such that $Z_N^{0,+,0}\leq
C^{\abs\thebox}$. To prove the latter inequality 
one can consider a shortest 
self-avoiding path $\go$ on $\bbZ^d$ starting at some
site of $\bnd\thebox$ and containing all the sites of $\thebox$, and using
$H^0_{0,N}(\phi) \geq \sum_{n=1}^{\abs\go-1}\Psi(\phi_{\go_n} - \phi_{\go_{n+1}})$.
%
%

\medskip

The above results imply the existence of
a wetting transition in these models\footnote{In the
$\gd$-pinning case, it can easily be seen that $\rho$ is monotonous in
$\gep$, so that there is a single critical value.}. 
Together with the
result of \cite{BDZ} that in the Gaussian model in
$d\geq 3$ there is {\em no} wetting transition, 
the first part of Theorem~\ref{thm_deltapinning} shows a radical difference
of behavior between the Gaussian and the SOS interactions.

\bigskip
\noindent
{\bf Remark:} 1. It is not difficult to see, looking at the proofs, that
our theorems remain true if we replace the SOS interaction $\Psi(x) = \abs x$
with any concave, even function, and the Gaussian interaction 
$\Psi(x)=\frac12 x^2$ with any even, convex 
$\Psi$ such that $1/c \geq \Psi''(x) \geq c$ for some $c>0$ and all $x$.\\
2. Even though the present work provides a proof of the wetting transition in
the models considered, several important issues remain completely open. In
particular, it would be most desirable to have a pathwise description of the
field in both the localized and repelled regimes, i.e. a proof that $\rho>0$
implies finiteness of the mean height of any fixed spin in the thermodynamic
limit (if possible with estimate on the tail and exponential decay of
correlations), and a proof that $\rho=0$ implies that the mean height of any
fixed spin diverges (if possible with estimates on the rate).
\subsection*{Acknowledgments} We thank Dima Ioffe and Ofer Zeitouni for their
comments and interest in this work.
\section{Proof of Theorem~\ref{thm_abs}}
Let $M>0$; following \cite{Chalker}, we introduce the set $\cB_M =
\setof{\phi}{\nu_N(\phi) \geq M}$, and the set
$$
\cC_M = \setof{\phi}{\sum_{x\in\thebox} \IF{\phi_x\leq 2a} \geq M}\,.
$$
Since $\cB_M \subset \cC_M$, the first claim immediately 
follows from the estimate on conditional probabilities
\begin{equation}
\mu_N^{0,+,a,b} (\cB_M | \cC_M) \leq e^{-C_2\,M}\,,\;\;\; M > C_1 N^{d-1}\,.
\label{conditional}
\end{equation}
Moreover, $\rho(\gep) = 0$ will follow from this and the obvious bound
\begin{equation*}\label{eq_obviousub}
\rho_N \leq \frac M {N^d} + \mu_N^{0,+,a,b}(\cB_M)\,,
\end{equation*}
by choosing $M$ such that $N^d\gg M\gg N^{d-1}$. 

\medskip
We turn to the proof of (\ref{conditional}).
We define a map $\frT$ from $\cC_M$ onto $\cB_M$ by
$$
(\frT\phi)_x = \begin{cases}
\phi_x   & \text{if }\phi_x \leq a\,,\\
\phi_x-a & \text{otherwise.}
\end{cases}
$$
If we write $e^{-V_N^{a,b}} = 
\IF{\phi_x>a} + e^b \IF{\phi_x\leq a}$ and expand the corresponding
products, we have
\begin{align*}
\mu&_N^{0,+,a,b}(\cC_M) = \int \frac{e^{-H_N^{0,a,b}(\phi)}}{Z_N^{0,+,a,b}}
\,\IF{\phi\in\cC_M}\,
\prod_{x\in\thebox}\dd\phi_x\\
& = \sumtwo{A\subset\thebox}{\abs A \geq M} \sum_{B\subset A} e^{b\abs B} \int
\frac{e^{-H_N^{0,a,b}(\phi)}}{Z_N^{0,+,a,b}}\prod_{x\in B} 
\IF{\phi_x\leq a} \dd\phi_x
\prod_{y\in A\setminus B} \IF{a<\phi_y\leq 2a} \dd\phi_y \prod_{z\in
\thebox\setminus A} \IF{2a<\phi_z}\dd\phi_z\,.
\end{align*}
Now, observe that
$$
H_{0,N}^0(\phi) \leq H_{0,N}^0(\frT\phi) + da\abs{\bnd\thebox} + 2da\abs{B}\,
$$
($\bnd\thebox$ being the set of $x\in\thebox$ neighbouring a site
$y\notin\thebox$). After the change of variables 
$\widetilde\phi_x=(\frT\phi)_x$, we have
\begin{align*}
\mu&_N^{0,+,a,b}(\cC_M)\\
&\geq  e^{-da\abs{\bnd\thebox}}\sumtwo{A\subset\thebox}{\abs A \geq
M} \sum_{B\subset A} \left(e^{-2da}e^b\right)^{\abs B} \int
\frac{e^{-H_{0,N}^0(\widetilde\phi)}}{Z_N^{0,+,a,b}} 
\prod_{x\in A} \IF{\widetilde\phi_x\leq a}
\dd\widetilde\phi_x\, \prod_{y\in\thebox\setminus A}
\IF{a<\widetilde\phi_y}\dd\widetilde\phi_y \nonumber\\
&= e^{-da\abs{\bnd\thebox}}\sumtwo{A\subset\thebox}{\abs A \geq
M} e^{b\abs A}\left(e^{-2da}+e^{-b}\right)^{\abs A} \int
\frac{e^{-H_{0,N}^0(\widetilde\phi)}}{Z_N^{0,+,a,b}}  
\prod_{x\in A} \IF{\widetilde\phi_x\leq a}
\dd\widetilde\phi_x\, \prod_{y\in\thebox\setminus A}
\IF{a<\widetilde\phi_y}\dd\widetilde\phi_y \nonumber\\
&\geq e^{-da\abs{\bnd\thebox}} \left(e^{-2da}+e^{-b}
\right)^M\;\mu_N^{0,+,a,b}(\cB_M)\,,
\label{eq_condabs}
\end{align*}
where we used $e^{-2da}+e^{-b}>1$, which follows from $ae^b < (2d)^{-1}$. 
This proves (\ref{conditional}). 
\section{Proof of Theorem~\ref{thm_deltapinning}}
The proof is very similar. Let $M>0$ and define $\cB_M$ as in the
previous proof (but remember that now $\nu_N$ is the number of sites with
height equal to $0$). We also need a set analogous to the set $\cC_M$ of the
previous section:
$$
\cD_M = \setof{\phi}{\sum_{x\in\thebox} \IF{\phi_x\leq 1} \geq M}\,.
$$
We are going to show that 
\begin{equation}
\mu_N^{0,+,\gep} (\cB_M | \cD_M) \leq e^{-C_4\,M}\,,\;\;\; M > C_3 N^{d-1}\,.
\label{conditional2}
\end{equation}
As in the previous theorem, (\ref{conditional2}) is sufficient 
to prove our claims.
\medskip

To prove (\ref{conditional2}) define the map
$$
(\frS\phi)_x = \begin{cases}
\phi_x-1 & \text{if }\phi_x>1\\
0        & \text{otherwise}
\end{cases}
$$
from $\cD_M$ onto $\cB_M$. 
Note that $\mu_N^{0,+,\gep}(\cD_M)$ can be written
\begin{equation}
\sumtwo{A\subset\thebox}{\abs A \geq M} \sum_{B\subset A} \gep^{\abs B} \int
\frac{e^{-H_{0,N}^0(\phi)}}{Z_N^{0,+,\gep}} 
\prod_{x\in A\setminus B} \IF{\phi_x\leq 1} \dd\phi_x
\prod_{y\in \thebox\setminus A} \IF{1<\phi_y}\dd\phi_y \prod_{z\in B}
\gd_0(\dd\phi_z)\,.
\label{eq_Zgauss}
\end{equation}
Let us first discuss the simpler case 
$\Psi(x)= \abs x$, $d \geq 1$. We have 
$$
H_{0,N}^0(\phi) \leq H_{0,N}^0(\frS\phi) + 
d\abs{\bnd\thebox} + 2d\abs{A}\,,
$$
and therefore, letting $\widetilde\phi_x=(\frS\phi)_x$ and
integrating over the variables $\phi_x$, $x\in A\setminus
B$,
\begin{align*}
\mu_N^{0,+,\gep}(\cD_M) &\geq e^{-d\abs{\bnd\thebox}}
\sumtwo{A\subset\thebox}{\abs A \geq
M}  \sum_{B\subset A} e^{-2d\abs A}\,\gep^{\abs B} \int
\frac{e^{-H_{0,N}^0(\widetilde\phi)}}{Z_N^{0,+,\gep}}
\prod_{x\in
\thebox\setminus A}\dd\widetilde\phi_x \prod_{y\in A}
\gd_0(\dd\phi_y)\\
&= e^{-d\abs{\bnd\thebox}}\sumtwo{A\subset\thebox}{\abs A \geq M}
\gep^{\abs A} e^{-2d\abs A}\,\left(1+\frac1\gep\right)^{\abs A} 
\frac{Z^{0,+}_\Ac}{Z_N^{0,+,\gep}}   \\
& \geq 
e^{-d\abs{\bnd\thebox}} (1+ e^{-2d})^M 
\,\mu_N^{0,+,\gep}(\cB_M)\,, 
\end{align*}
where we have used the assumption
$\gep < e^{-2d}$. This proves (\ref{conditional2}).

\medskip

Let us now turn to the case $\Psi(x) = \frac12 x^2$, $d = 2$.
Writing $W=A\cup\thebox^c$, we have the estimate
$$
H_{0,N}^0(\phi) \leq H_{0,N}^0(\frS\phi) + 2\abs{\bnd\thebox} + 
8\abs{A} + 2 \sum_{x\in\bnd W}\sumtwo{y\notin W}{y\nn x} (\frS\phi)_y\,.
$$
Let us use the short-hand notation 
$\cX(\phi) = 2 \sum_{x\in\bnd W}\sum_{y\notin W,\,y\nn x} \phi_y$.
Inserting this estimate in
\eqref{eq_Zgauss} 
and changing variables to $\widetilde\phi_x = (\frS\phi)_x$, we obtain
\begin{equation*}
\mu_N^{0,+,\gep}(\cD_M) 
\geq e^{-2\abs{\bnd\thebox}}
\sumtwo{A\subset\thebox}{\abs A \geq M}
\gep^{\abs A}\; (e^{-8}(1+\frac1\gep))^{\abs A}\;  
\frac{Z^{0,+}_\Ac}{Z_N^{0,+,\gep}}\;
\mu^{0,+}_\Ac \left(e^{-\cX(\phi)}\right)\,.
\end{equation*}
Clearly,
$$
\mu^{0,+}_\Ac \left(e^{-\cX(\phi)}\right) \geq
e^{-c_1(\abs{\bnd\thebox} + \abs A)}\; \mu^{0,+}_\Ac(\cX(\phi) \leq c_1\abs{\bnd W})\,.
$$
The conclusion will follow as before, once we prove that this last probability
is bounded from below by, say, $1/2$. Using Markov inequality, we get
$$
\mu^{0,+}_\Ac(\abs{\bnd W}^{-1}\cX(\phi) >
c_1) \leq \frac{\mu^{0,+}_\Ac(\abs{\bnd W}^{-1}\cX(\phi))}{c_1}\,.
$$
Since it follows from~\cite{DMRR} (Lemma 2.1) that there exists an absolute
constant $c_2$ such that $\mu^{0,+}_\Ac(\abs{\bnd W}^{-1}\cX(\phi)) \leq c_2$,
the result follows by taking $c_1$ large enough.
\end{document}